\documentclass[12pt]{amsart}

\newtheorem{theorem}{Theorem}
\newtheorem{proposition}[theorem]{Proposition}
\newtheorem{corollary}[theorem]{Corollary}

\theoremstyle{definition}
\newtheorem{definition}[theorem]{Definition}

\theoremstyle{remark}
\newtheorem{remark}[theorem]{Remark}

\newcommand{\tensor}{\otimes}

\newcommand{\CC}{\mathbb{C}}
\newcommand{\PP}{\mathbb{P}}
\newcommand{\RR}{\mathbb{R}}

\newcommand{\C}{\mathcal{C}}
\newcommand{\E}{\mathcal{E}}
\newcommand{\F}{\mathcal{F}}

\renewcommand{\O}{\mathcal{O}}

\DeclareMathOperator{\Ext}{Ext}
\DeclareMathOperator{\Hom}{Hom}
\DeclareMathOperator{\id}{id}
\DeclareMathOperator{\rank}{rank}
\DeclareMathOperator{\SL}{SL}
\DeclareMathOperator{\SU}{SU}
\DeclareMathOperator{\tr}{tr}
\DeclareMathOperator{\U}{U}
\DeclareMathOperator{\vol}{vol}

\begin{document}

\baselineskip=15pt

\title[Vortex equation and reflexive sheaves]{Vortex
equation and reflexive sheaves}

\author[I. Biswas]{Indranil Biswas}

\address{School of Mathematics, Tata Institute of Fundamental
Research, Homi Bhabha Road, Bombay 400005, India}

\email{indranil@math.tifr.res.in}

\author[M. Stemmler]{Matthias Stemmler}

\email{stemmler@math.tifr.res.in}



\date{}

\begin{abstract}
It is known that given a stable holomorphic pair $(E\, ,\phi)$, where $E$ is a holomorphic vector bundle on a compact K\"ahler manifold $X$ and $\phi$ is a holomorphic section of $E$, the vector bundle $E$ admits a Hermitian metric solving the vortex equation. We generalize this to pairs $(\E\, ,\phi)$, where $\E$ is a reflexive sheaf on $X$.
\end{abstract}

\maketitle

\section{Introduction}

The Hitchin-Kobayashi correspondence states that a holomorphic vector bundle $E$ on a compact K\"ahler manifold $(X\, , \omega)$ admits a Hermitian-Einstein metric if and only if it is polystable in the sense of Mumford-Takemoto. In \cite{Br90}, \cite{Br91}, \cite{GP93} and \cite{GP94a}, Bradlow and Garc\'ia-Prada established a generalization of this correspondence to the context of holomorphic pairs, consisting of a holomorphic vector bundle $E$ and a global holomorphic section $\phi$ which is not identically equal to zero. Instead of Hermitian-Einstein metrics, they considered Hermitian metrics on $E$ satisfying the so-called $\tau$-vortex equation, which additionally depends on the section $\phi$ and a real parameter $\tau$. This equation can be seen as a dimensional reduction of the Hermitian-Einstein equation for an $\SU(2)$-equivariant vector bundle on $X \times \PP^1$, where $\SU(2)$ acts trivially on $X$ and in the standard way on $\PP^1$ \cite{GP94b}. The notion of a stable vector bundle on $X$ extends to the notion of $\tau$-stable holomorphic pairs. Using the dimensional reduction procedure it was shown that a $\tau$-stable holomorphic pair admits a Hermitian metric satisfying the $\tau$-vortex equation.

In \cite{BS94}, Bando and Siu generalized the notion of Hermitian-Einstein metrics to reflexive sheaves $\E$ by considering a special class of Hermitian metrics on the locally free part of $\E$ called admissible metrics. They proved that every stable reflexive sheaf on a compact K\"ahler manifold admits an admissible Hermitian-Einstein metric. This work of Bando and Siu has turned out to be extremely useful testified by its numerous applications.

Our aim here is to unify these two generalizations to obtain an analogous existence theorem for pairs consisting of a reflexive sheaf 
and a global section. This answers a question of Tian and Yang (see Remark 1 in \cite{TY02} after \cite[Theorem 4.6]{TY02}).

We prove the following (see Theorem \ref{result} and Corollary \ref{cor1}):

\begin{theorem}
Let $(\E\, , \phi)$ be a reflexive sheaf pair on a compact K\"ahler manifold $(X\, , \omega)$, and let $\tau$ be a real number. Let $S 
\subset X$ be the singular set for $\E$, and
\[
  \widehat{\tau} \,=\, \frac{\tau \vol(X)}{4 \pi}\, .
\]
Then $(\E\, , \phi)$ is $\widehat \tau$-polystable if and only if there exists an admissible Hermitian metric on $\E|_{X \setminus S}$ 
satisfying the $\tau$-vortex equation.
\end{theorem}

In \cite{BS} a corresponding result for Higgs sheaves was established.

\section{Preparations and statement of the theorem}

Let $(X\, , \omega)$ be a compact connected K\"ahler manifold of complex dimension $n$, and let $\E$ be a torsion-free coherent analytic sheaf on $X$. We recall a few standard definitions.
\begin{definition} \mbox{}
\begin{enumerate}
\item[(1)] The {\em degree\/} of $\E$ (with respect to $\omega$) is defined in terms of cohomology classes as follows:
\[
  \deg(\E) \,:=\, \frac{1}{(n-1)!} \, (c_1(\E) \cup [\omega]^{n-1}) \cap [X] \in \RR\, ,
\]
where $c_1(\E) = c_1(\det \E)$ is the first Chern class of $\E$, defined using the determinant line bundle $\det \E$ of $\E$, and $[\omega] \in H^2(X, \RR)$ is the cohomology class of $\omega$.

\item[(2)] The {\em singular set\/} for $E$ is the closed analytic subset of $X$ outside which $E$ is locally free.

\item[(3)] If $\rank(\E) > 0$, the {\em slope\/} of $\E$ is defined to be the ratio
\[
  \mu(\E) \, :=\, \frac{\deg(\E)}{\rank(\E)}\, .
\]
Here, $\rank(\E)$ denotes the rank of $\E$ outside its singular set.
\end{enumerate}
\end{definition}

\begin{definition} \mbox{}
\begin{enumerate}
\item[(1)] $\E$ is said to be {\em stable\/} (with respect to $\omega$) if
\[
  \mu(\E') < \mu(\E)
\]
holds for every coherent analytic
subsheaf $\E'$ of $\E$ with $0 < \rank(\E') < \rank(\E)$.
\item[(2)] $\E$ is said to be {\em polystable\/} if $\E$ is a direct sum
\[
  \E = \E_1 \oplus \cdots \oplus \E_m
\]
of stable subsheaves with the same slope $\mu(\E_i) = \mu(\E)$ for all $i$.
\end{enumerate}
\end{definition}

For later use, we also explain how this relates to group actions. Let $G$ be a compact Lie group acting holomorphically on $X$ and preserving the K\"ahler form $\omega$. Then there is a $G$-invariant version of the above definition. Namely, let $\E$ be $G$-equivariant, meaning $\E$ is equipped with a lift of the action of $G$ on $X$ to $\E$.

\begin{definition}
$\E$ is said to be {\em $G$-invariantly stable\/} (with respect to $\omega$) if
\[
  \mu(\E') < \mu(\E)
\]
holds for every $G$-invariant coherent analytic
subsheaf $\E'$ of $\E$ with $0 < \rank(\E') < \rank(\E)$.
\end{definition}

Now let $E$ be a holomorphic vector bundle on $X$.
\begin{definition}
A Hermitian metric $h$ on $E$ is called a {\em Hermitian-Einstein metric\/} (with respect to $\omega$) if
\[
  \sqrt{-1} \Lambda_\omega F_h = \lambda \id_E \quad \text{for some } \lambda \in \RR\, ,
\]
where $\Lambda_\omega$ is the adjoint of forming the wedge product with $\omega$, $F_h$ is the curvature of the Chern connection for $h$ and $\id_E$ is the identity automorphism of $E$.
\end{definition}

There is a correspondence between stability and Hermitian-Einstein metrics known as the Hitchin-Kobayashi correspondence, which was proved by Donaldson (\cite{Do85}, \cite{Do87}) and Uhlenbeck and Yau (\cite{UY86}, \cite{UY89}). We recall it:

\begin{theorem}
A holomorphic vector bundle $E$ on $X$ admits a Hermitian-Einstein metric if and only if $E$ is polystable.
\end{theorem}

In \cite{BS94}, Bando and Siu established a Hitchin-Kobayashi correspondence for the more general situation of reflexive sheaves.

Let $\E$ be a torsion-free coherent analytic
sheaf on $X$. Let $S \,\subset\,
X$ be the singular set for $\E$. So $S$ is a closed complex
analytic subset of $X$ of codimension at least $2$. We recall
the definition of Hermitian-Einstein metrics for sheaves.

\begin{definition} \label{def-ad}
A Hermitian metric $h$ on the holomorphic vector bundle $\E|_{X \setminus S}$ is called {\em admissible\/} if the following conditions are satisfied:
\begin{enumerate}
\item[(A1)] The Chern curvature $F_h$ of $h$ is square-integrable.
\item[(A2)] The contracted Chern curvature $\Lambda_\omega F_h$ is bounded.
\end{enumerate}
\end{definition}

Using a heat equation approach, Bando and Siu proved the following version of the Hitchin-Kobayashi correspondence.
\begin{theorem}[\cite{BS94}, Theorem 3] \label{Bando-Siu}
A reflexive sheaf $\E$ on a compact K\"ahler manifold admits an admissible Hermitian-Einstein metric if and only if $\E$ is polystable.
\end{theorem}

In what follows, we will be concerned with what we call sheaf pairs.
\begin{definition}
A {\em sheaf pair\/} on $X$ is a pair $(\E\, , \phi)$, where
\begin{itemize}
\item $\E$ is a coherent analytic sheaf on $X$, and
\item $\phi \in H^0(X, \E) \setminus \{ 0 \}$ is a global section of $\E$ (so $\phi$ is not identically zero).
\end{itemize}
A sheaf pair $(\E\, , \phi)$ is called {\em torsion-free\/} (respectively, {\em reflexive}) if $\E$ is a torsion-free (respectively, reflexive) sheaf on $X$. If $\E$ is locally free, then $(\E\, , \phi)$ is called a {\em holomorphic pair}.
\end{definition}

In \cite{Br91}, Bradlow introduced a notion of stability for
holomorphic pairs (see also \cite{GP94b}) depending on a real
parameter $\tau$; the definition extends to the 
set-up of torsion-free sheaf pairs.

\begin{definition}
Let $\tau$ be a real number. A torsion-free sheaf pair $(\E\, , \phi)$ is called {\em $\tau$-stable\/} (with respect to $\omega$) if the following conditions are satisfied.
\begin{itemize}
\item $\mu(\E') < \tau$ for every coherent analytic
subsheaf $\E'$ of $\E$ with $\rank(\E') > 0$.
\item $\mu(\E/\E') > \tau$ for every coherent analytic
subsheaf $\E'$ of $\E$ with $0 < \rank(\E') < \rank(\E)$ and $\phi \in H^0(X, \E')$.
\end{itemize}
A torsion-free sheaf pair $(\E\, , \phi)$ is called {\em $\tau$-polystable\/} if it is either $\tau$-stable or $\E$ decomposes as a direct sum
\[
  \E \,=\, \E' \oplus \E''
\]
of coherent analytic subsheaves such that $\phi \in H^0(X, \E')$, the sheaf pair $(\E'\, , \phi)$ is $\tau$-stable and the sheaf $\E''$ is polystable with slope $\mu(\E'') = \tau$.
\end{definition}

In the case of holomorphic pairs, the appropriate replacement for a Hermitian-Einstein metric is a metric satisfying the so-called vortex equation, which also depends on a real parameter $\tau$ and has been studied by Bradlow \cite{Br90}, \cite{Br91} and Garc\'ia-Prada \cite{GP93}, \cite{GP94a}.

\begin{definition}
Let $\tau$ be a real number. Given a holomorphic pair $(E\, , \phi)$, a Hermitian metric $h$ on $E$ is said to satisfy the {\em $\tau$-vortex equation\/} if
\[
  \Lambda_\omega F_h - \frac{\sqrt{-1}}{2} \, \phi \circ \phi^\ast + \frac{\sqrt{-1}}{2} \, \tau \id_E = 0\, ,
\]
where $\phi$ is regarded as a homomorphism from the trivial holomorphic Hermitian line bundle on $X$ to $E$ and $\phi^\ast$ denotes the adjoint of $\phi$ with respect to $h$.
\end{definition}

The following theorem is proved in \cite{GP94b} (see \cite[Theorem 4.33]{GP94b}).

\begin{theorem} \label{Garcia-Prada}
Let $(E\, , \phi)$ be a holomorphic pair on a compact K\"ahler manifold $(X\, , \omega)$, and let $\tau$ be a real number. If $(E\, , \phi)$ is $\widehat \tau$-stable, where
\[
  \widehat{\tau} = \frac{\tau \vol(X)}{4 \pi}\, ,
\]
then there exists a Hermitian metric on $E$ satisfying the $\tau$-vortex equation. Here, $\vol(X)$ denotes the volume of $X$ with respect to $\omega$.
\end{theorem}

We will prove the following theorem.

\begin{theorem} \label{result}
Let $(\E\, , \phi)$ be a reflexive sheaf pair on a compact K\"ahler manifold $(X\, , \omega)$, and let $\tau$ be a real number. Let $S \subset X$ be the singular set for $\E$, and define
\[
  \widehat{\tau} \,=\, \frac{\tau \vol(X)}{4 \pi}\, .
\]
If $(\E\, , \phi)$ is $\widehat\tau$-stable, then there exists an admissible Hermitian metric on $\E|_{X \setminus S}$ satisfying the $\tau$-vortex equation.
\end{theorem}

\begin{remark}
If $\phi\,=\,0$, then the $\tau$-vortex equation imposes the constraint that $\widehat \tau \,=\,\mu(\E)$. However, given a stable reflexive sheaf $\E$, the pair $(\E\, , 0)$ is $\tau'$-stable for all $\tau' = \mu(\E) + \varepsilon$ with $\varepsilon > 0$ sufficiently small, implying that the above constraint is violated. This accounts for the condition $\phi \neq 0$ in the definition of sheaf pairs.
\end{remark}

The rest of this note will be concerned with the proof of theorem \ref{result}. We first apply the technique of dimensional reduction to obtain an $\SU(2)$-invariantly stable reflexive sheaf on $X \times \PP^1$. Then we show that the admissible Hermitian-Einstein metric obtained by Theorem \ref{Bando-Siu} yields an admissible Hermitian metric satisfying the vortex equation.

\section{Proof of Theorem \ref{result}}

Suppose we are in the situation of Theorem \ref{result}. First we apply the technique of dimensional reduction for sheaf pairs, which was developed in \cite{GP94b}.

Consider the compact complex manifold $X \times \PP^1$, where $\PP^1$ is the complex projective line. Let
\[
  p\, :\, X \times \PP^1\,\longrightarrow\, X ~\,~\text{~and~}~\,~
  q\, :\, X \times \PP^1\,\longrightarrow\, \PP^1
\]
be the natural projections. Let $\SU(2)$ act on $X \times \PP^1$ by
acting trivially on $X$ and in the standard way on $\PP^1$, by
regarding $\PP^1$ as $({\mathbb C}^2 \setminus\{0\})/{\mathbb C}^*
\, =\, \SU(2)/\U(1)$.

According to \cite{GP94b} (see \cite[proof of Theorem 4.9]{GP94b}), the sheaf pair $(\E\, , \phi)$ defines an
$\SU(2)$-equivariant coherent analytic sheaf
$\F$ on $X \times \PP^1$ given as an extension
\begin{equation} \label{extension}
  0 \longrightarrow p^\ast \E  \longrightarrow  \F \longrightarrow q^\ast \O(2) \longrightarrow 0\, ,
\end{equation}
where $\O(2)\,=\, T\PP^1$ is the line bundle of degree $2$ on $\PP^1$ (see also \cite[Theorem 1.1]{AG01}). To see how the section $\phi$ is involved, note that extensions of the form \eqref{extension} are parametrized by
\[
  \Ext_{X \times \PP^1}^1(q^\ast \O(2), p^\ast \E)\, .
\]
Since $q^\ast \O(2)$ is locally free, this group is isomorphic to $H^1(X \times \PP^1, p^\ast \E \tensor q^\ast \O(-2))$, which is isomorphic to $H^0(X, \E)$ (cf.\ \cite[proof of Theorem 4.9]{GP94b}).

Since $\E$ is in particular a torsion-free sheaf on $X$, the complex analytic subset $S$ in Theorem \ref{result} is of codimension at least $2$. Being an extension of $p^\ast \E$ by the locally free sheaf $q^\ast \O(2)$, the sheaf $\F$ is locally isomorphic to the direct sum $p^\ast \E \oplus q^\ast \O(2)$. Therefore, $\F$ is locally free on $p^{-1}(X \setminus S) \,=\, (X \setminus S) \times \PP^1$. Also, since the projection $p$ is a flat morphism, $\F$ is a reflexive sheaf on $X \times \PP^1$.

Denote by $E$ and $F$ the holomorphic vector bundles corresponding to the locally free sheaves $\E|_{X \setminus S}$ and $\F|_{(X \setminus S) \times \PP^1}$, respectively. Note that in the $\C^\infty$ category, we have an isomorphism
\[
  F \simeq p^\ast E \oplus q^\ast \O(2)\, .
\]

There is a relation between the $\tau$-stability of the pair $(\E\, , \phi)$ and the stability of $\F$ with respect to a special K\"ahler metric on $X \times \PP^1$ encoding the parameter $\tau$. More precisely, for every positive real number $\sigma$ consider the K\"ahler form 
\[
  \Omega_\sigma \,:=\, p^\ast \omega \oplus \sigma q^\ast \omega_{\PP^1}
\]
on $X \times \PP^1$, where $\omega_{\PP^1}$ is the normalized Fubini-Study K\"ahler form on $\PP^1$, such that $\int_{\PP^1} \omega_{\PP^1} = 1$. Then we have the following theorem \cite[Theorem 4.21]{AG03} (see also \cite[Theorem 4.9]{GP94b}):

\begin{theorem}
Let $(\E\, , \phi)$ be a sheaf pair on a compact K\"ahler manifold $(X\, , \omega)$. Let $\F$ be
the $\SU(2)$-equivariant coherent
analytic sheaf on $X \times \PP^1$ determined by $(\E\, , \phi)$ as the extension \eqref{extension}, and let the numbers $\sigma$ and $\tau$ be related by
\begin{equation} \label{sigma-tau}
  \sigma = \frac{2 \vol(X)}{(\rank(\E) + 1) \, \tau - \deg(\E)}\, .
\end{equation}
Then $(\E\, , \phi)$ is $\tau$-stable if and only if $\sigma > 0$ and $\F$ is $\SU(2)$-invariantly stable with respect to $\Omega_\sigma$.
\end{theorem}

Consequently, in the situation of Theorem \ref{result}, we know that $\F$ is $\SU(2)$-invariantly stable with respect to $\Omega_\sigma$, where $\sigma$ is determined from $\tau$ by \eqref{sigma-tau}. As in \cite[Theorem 6]{GP93}, it follows that $\F$ is polystable with respect to $\Omega_\sigma$. By Theorem \ref{Bando-Siu}, there is an admissible Hermitian-Einstein metric $\widetilde h$ on the holomorphic vector bundle $F$ with respect to $\Omega_\sigma$. By pulling back $\widetilde h$ by each element of $\SU(2)$ and averaging over the group using the Haar measure on the compact group $\SU(2)$, we can assume that $\widetilde h$ is an $\SU(2)$-invariant Hermitian-Einstein metric, cf.\ \cite[proof of Theorem 5]{GP93}. As in \cite[Proposition 3.2]{GP94b}, according to the $\C^\infty$ decomposition
\begin{equation} \label{decomposition}
  F \simeq p^\ast E \oplus q^\ast \O(2)
\end{equation}
of vector bundles on $(X \setminus S) \times \PP^1$, the metric $\widetilde h$ decomposes as
\[
  \widetilde{h}\, =\, p^\ast h \oplus q^\ast h'\, ,
\]
where $h$ is a Hermitian metric on the holomorphic vector bundle $E$ on $X \setminus S$, and $h'$ is an $\SU(2)$-invariant Hermitian metric on the holomorphic line bundle $\O(2)$ on $\PP^1$. In order to complete the proof of the theorem, we have to show that $h$ is admissible and satisfies the $\tau$-vortex equation. The latter follows as in \cite[Proposition 3.11]{GP94b}.

\begin{proposition}
Let the numbers $\sigma$ and $\widehat \tau$ be related by
\[
  \sigma = \frac{2 \vol(X)}{(\rank(E) + 1) \, \widehat{\tau} - \deg E}\, ,
\]
where
\[
  \widehat{\tau} = \frac{\tau \vol(X)}{4 \pi}\, .
\]
Then in the situation considered above, the metric $h$ on $E$ satisfies the $\tau$-vortex equation if and only if $\widetilde h$ is a Hermitian-Einstein metric on $F$ with respect to $\Omega_\sigma$.
\end{proposition}

We repeat some details of the argument so that we can explain how it also yields the admissibility of $h$. By \cite[Proposition 3.5]{GP94b}, the Chern connection $D_{\widetilde h}$ for the Hermitian holomorphic vector bundle $(F\, , \widetilde{h})$ can be written according to the decomposition \eqref{decomposition} as
\[
  D_{\widetilde h} = \begin{pmatrix}
    D_{p^\ast h} & \beta         \\
    - \beta^\ast & D_{q^\ast h'}
  \end{pmatrix},
\]
where $D_{p^\ast h}$ and $D_{q^\ast h'}$ are the Chern connections for $p^\ast h$ and $q^\ast h'$, respectively, $\beta$ is a $1$-form with values in $\Hom(q^\ast \O(2), p^\ast E)$, and $\beta^\ast$ is the adjoint of $\beta$ with respect to $\widetilde h$. In fact,
\[
  \beta = p^\ast \phi \tensor q^\ast \alpha\, ,
\]
where $\alpha \in \Omega^{0,1}(\PP^1, \O(-2))$ is an $\SU(2)$-invariant form which is unique up to a multiplicative constant. The curvature of $D_{\widetilde h}$ can then be written as
\begin{equation} \label{curvature decomposition}
  F_{\widetilde h} = \begin{pmatrix}
    F_{p^\ast h} - \beta \wedge \beta^\ast & \partial \beta                          \\
    - \bar \partial \beta^\ast             & F_{q^\ast h'} - \beta^\ast \wedge \beta
  \end{pmatrix},
\end{equation}
where $\partial$ and $\bar \partial$ denote the components of the induced connection on $\Hom(q^\ast \O(2), p^\ast E)$. If $\alpha$ is chosen such that $\alpha \wedge \alpha^\ast = \frac{\sqrt{-1}}{2} \, \sigma \omega_{\PP^1}$, then we have
\begin{align*}
  \beta \wedge \beta^\ast &= \frac{\sqrt{-1}}{2} \, \sigma p^\ast(\phi \circ \phi^\ast) \tensor q^\ast \omega_{\PP^1}\, ,   \\
  \beta^\ast \wedge \beta &= - \frac{\sqrt{-1}}{2} \, \sigma p^\ast(\phi^\ast \circ \phi) \tensor q^\ast \omega_{\PP^1}\, .
\end{align*}
The Hermitian-Einstein equation for $\widetilde h$,
\[
  \sqrt{-1} \Lambda_{\Omega_\sigma} F_{\widetilde h} = \lambda \id_F\, ,
\]
can then be translated into the equations
\begin{align}
  \label{vortex 1} \Lambda_\omega F_h - \frac{\sqrt{-1}}{2} \, \phi \circ \phi^\ast + \frac{\sqrt{-1}}{2} \, \tau \id_E    &= 0\, , \\
  \label{vortex 2} \frac{\sqrt{-1}}{2} \, \phi^\ast \circ \phi - \frac{4 \pi \sqrt{-1}}{\sigma} + \frac{\sqrt{-1}}{2} \, \tau &= 0\, ,
\end{align}
and \eqref{vortex 1} is the vortex equation for $(E\, , h)$.

We shall check the admissibility of $h$. Condition (A2) in Definition \ref{def-ad} follows immediately from \eqref{vortex 1} and \eqref{vortex 2}. In fact, \eqref{vortex 2} implies that the function
\[
  \phi^\ast \circ \phi = \tr(\phi^\ast \circ \phi) \,=\, \tr(\phi \circ \phi^\ast) \,=\, |\phi|_h^2
\]
is constant, where $\tr$ denotes the trace and $\phi^\ast \circ \phi$ and $\phi \circ \phi^\ast$ are regarded as endomorphisms of $\O_X$ and $E$ respectively. Consequently, the norm of $\phi \circ \phi^\ast$ is
\[
  |\phi \circ \phi^\ast|_h^2 \,=\, \tr(\phi \circ \phi^\ast \circ (\phi \circ \phi^\ast)^\ast) \,= \,|\phi|_h^4
\]
and therefore it is also constant. Now \eqref{vortex 1} implies the boundedness of $|\Lambda_\omega F_h|_h$.

Condition (A1) follows in the same way from the admissibility of $\widetilde h$: Since $F_{\widetilde h}$ is square-integrable, we know by \eqref{curvature decomposition} that $F_{p^\ast h} - \beta \wedge \beta^\ast$ and $F_{q^\ast h'} - \beta^\ast \wedge \beta$ are square-integrable. Since $F_{q^\ast h'}$ has no singularities and the norms of $\beta \wedge \beta^\ast$ and $\beta^\ast \wedge \beta$ coincide, it follows that $F_{p^\ast h}$, and hence $F_h$, is square-integrable. This shows the admissibility of $h$, and completes the proof of Theorem \ref{result}.

It is straight-forward to check that if a reflexive sheaf pair $(\E\, , \phi)$ on $(X\, , \omega)$ has an admissible Hermitian metric satisfying the $\tau$-vortex equation, then $(\E\, , \phi)$ is $\widehat\tau$-polystable. Therefore, Theorem \ref{result} has the following corollary.

\begin{corollary} \label{cor1}
Let $(\E\, , \phi)$ be a reflexive sheaf pair on a compact K\"ahler manifold $(X\, , \omega)$, and let $\tau$ be a real number. Let $S \subset X$ be the singular set for $\E$, and
\[
  \widehat{\tau} \,=\, \frac{\tau \vol(X)}{4 \pi}\, .
\]
Then $(\E\, , \phi)$ is $\widehat \tau$-polystable if and only if there exists an admissible Hermitian metric on $\E|_{X \setminus S}$ satisfying the $\tau$-vortex equation.
\end{corollary}

\end{document}